\documentclass[preprint,authoryear,12pt]{elsarticle}

\usepackage{amssymb,amsthm}
\usepackage[latin1]{inputenc}

\newtheorem{thm}{Theorem}
\newproof{pf}{Proof}

\journal{Computer Aided Geometric Design}

\begin{document}

\begin{frontmatter}



\title{The inversion problem for rational Bézier curves}


\author{Ana Marco\corref{1}}
\ead{ana.marco@uah.es}

\author{José-Javier Mart{\'\i}nez \corref{2}}
\ead{jjavier.martinez@uah.es}

\cortext[2]{Corresponding author}

\address{Departamento de Matemáticas, Universidad de Alcalá, Campus Universitario, 28871 Alcalá de Henares, Madrid, Spain}

\begin{abstract}

The inversion problem for rational Bézier curves is addressed by using resultant matrices for polynomials expressed in the Bernstein basis. The aim of the work is not to construct an inversion formula but finding the corresponding value of the parameter for each point on the curve. Since sometimes one has only an approximation of that point the use of the singular value decomposition, a key tool in numerical linear algebra, is shown to be adequate.
\end{abstract}


\end{frontmatter}


\section{Introduction}
\label{}

The inversion problem (given a point $P_0=(x_0,y_0)$ on a curve $C$ parametrized by $P(t)=(x(t),y(t))$), compute the value $t_0$ of the parameter $t$ corresponding to this point $P_0$) is a classical problem in the area of Computer Aided Geometric Design (CAGD) which has applications, for instance, in computing the intersection points of two curves \citep{HL}.

Two pioneering papers on this subject, which also study the problem of implicitization, are \citep{SAG, GSA}, the use of resultant matrices being a fundamental tool in those papers. A nice recent survey of those and related problems can be found in \citep{SZ}.

Two recent alternative approaches to the inversion problem have been presented in \citep{GVR,BDA}.

Our aim in this paper is to solve numerically the inversion problem for the class of Bézier curves, a very popular type of curves in CAGD \citep{FA}. Since for these curves the parametric equations are expressed in the Bernstein basis, it is convenient to construct the resultant matrices using the Bernstein basis directly, an idea already suggested in \citep{GSA}. 

A resultant matrix for polynomials expressed in the Bernstein basis is the so-called Bernstein-Bézout matrix, which has been used in \citep{BG} for designing a fast algorithm for computing the greatest common divisor of two polynomials expressed in the Bernstein basis, and then in \citep{MM07} for computing the implicit equation of a Bézier curve by means of interpolation. An algorithm for the fast computation of the Bézout resultant matrix for polynomials expressed in the monomial basis is introduced in \citep{CZG}, and an algorithm for the fast computation of the Bernstein-Bézout resultant matrix (which we will use in our approach to the inversion problem) is presented in \citep{BG}.

Another fundamental tool in our approach is the singular value decomposition (SVD). The SVD is one of the most valuable tools in numerical linear algebra and we will employ it in this work for computing an approximation of a basis of the nullspace of a Bernstein-Bézout matrix. The use of the SVD is specially appropriate when the matrix is not know exactly, what happens in our case when only an approximation of the point $P_0=(x_0,y_0)$ is known. It must be observed that this is the usual situation, for example, when the inversion problem is utilized to solve problems of curve intersection \citep{MM04}.

Related to this situation (the fact that the point is usually only close to, but not exactly on, the curve) in Section 10 of \citep{GSA} the problem of the geometric interpretation of {\it inversion formulae} was briefly addressed, and this problem was later deeply studied in \citep{WJ}.

It is important to point out that our approach does not carry out any transformation between Bernstein and power (monomial) basis, because this conversion could involve a greater loss of accuracy \citep{F}. In this sense, in \citep{BG} is indicated that {\it for numerical computations involving polynomials in Bernstein form it is essential to consider algorithms which express all intermediate results using this form only}.

The rest of the paper is organized as follows. In Section 2 some basic results on the Bernstein-Bézout matrix and the SVD are presented. In Section 3 our procedure for solving the inversion problem for a Bézier curve is introduced, leaving for Section 4 its analysis in the more general case of rational curves expressed in the Bernstein basis. Section 5 contains some numerical experiments that illustrate the performance of our approach, and finally Section 6 is devoted to conclusions.

\section{Some basic results}
\label{}
In this section we recall some basic results on two fundamental tools that we will use in this work: the Bernstein-Bézout matrix and the singular value decomposition.

\subsection{The Bernstein-Bézout matrix}

The Bernstein-Bézout matrix is a resultant matrix for polynomials expressed in the Bernstein basis. It can be defined in an analogous way to the definition of the Bézout resultant matrix for polynomials expressed in the monomial basis \citep{SZ}. 

Let $p(t)=p_0 \beta_0^{(n)}(t)+p_1 \beta_1^{(n)}(t)+ \cdots + p_n \beta_n^{(n)}(t)$ y $q(t)=q_0 \beta_0^{(n)}(t)+q_1 \beta_1^{(n)}(t)+ \cdots + q_n \beta_n^{(n)}(t)$ be two polynomials expressed in the Bernstein basis
$$
\mathcal B_n =\big\{ \beta_i^{(n)}(t) = {n \choose i} (1 - t)^{n-i} t^i,
\qquad i = 0, \ldots, n \big\}.
$$
The {\it Bernstein-Bézout matrix} $B=(b_{i,j}) \in R^{n\times n}$ of $p(t)$ and $q(t)$ is defined by
$$
\frac{p(t)q(s)-p(s)q(t)}{t-s}=\sum_{i,j=1}^n b_{i,j} \beta_{i-1}^{(n-1)}(t)\beta_{j-1}^{(n-1)}(s),
$$
which can be equivalently rewritten as
$$
\frac{p(t)q(s)-p(s)q(t)}{t-s}=\left(
                                \begin{array}{ccc}
                                  \beta_0^{(n-1)}(s) & \cdots & \beta_{n-1}^{(n-1)}(s) \\
                                \end{array}
                              \right)B\left(
                                        \begin{array}{c}
                                          \beta_0^{(n-1)}(t) \\
                                          \vdots \\
                                          \beta_{n-1}^{(n-1)}(t) \\
                                        \end{array}
                                      \right).
$$

The following result is a direct consequence of the definition of the Bernstein-Bézout matrix:

\begin{thm}
Let $p(t)$ and $q(t)$ be two polynomials expressed in the Bernstein basis $\mathcal B_n$, and $B$ be the Bernstein-Bézout matrix of $p(t)$ and $q(t)$. If $t_0$ is a common root of $p(t)$ and $q(t)$, then $det(B)=0$.
\end{thm}

A fast algorithm for computing the entries of the Bernstein-Bézout matrix (an algorithm which we will use in this work) has been presented in \citep{BG}. Since, for reasons explained in Section 5, we will use for the examples of that section the symbolic computation system {\it Maple} for the computation of the Bernstein-Bézout matrix, we include here the the Bini and Gemignani algorithm written in the {\it Maple} language:
\begin{verbatim}
for i from 1 to n do
   B[i,1]:=(n/i)*(p[i]*q[0]-p[0]*q[i]);
od;

for j from 1 to n-1 do
   B[n,j+1]:= (n/(n-j))*(p[n]*q[j]-p[j]*q[n])
od;

for j from 1 to n-1 do
   for i from 1 to n-1 do
     B[i,j+1]:=(n^2/(i*(n-j)))*(p[i]*q[j]-p[j]*q[i])
     +((j*(n-i))/(i*(n-j)))*B[i+1,j];
   od;
od;
\end{verbatim}

\subsection{The Singular Value Decomposition}

Let us now recall the concept of singular value decomposition
(SVD):

Given an $m \times n$ real matrix $A$,
there exist orthogonal matrices $U$ (of order $m$) and $V$ (of
order $n$), and a diagonal matrix $\Sigma$ (of size $m \times n$)
such that
$$
A = U \Sigma V^{T}.
$$
This factorization of $A$ is called the {\it singular value
decomposition (SVD)} of $A$.

The  $r$ (with $r \le m,n$) nonzero diagonal entries of $\Sigma$
are the {\it singular values} of $A$ ( i.e. the positive
square roots of the eigenvalues of $A^{T}A$). If there are $r$
(nonzero) singular values, then $r$ is the {\it rank} of $A$.

But we also obtain a very important additional advantage from the
computation of the SVD of a matrix $A$ of size $m \times n$
(including $V$, not only $\Sigma$): {\it the last $n-r$ columns of
of $V$ form a basis of the nullspace of $A$} (see \citep{ST, D}).

The SVD provides the best way of estimating the rank of a matrix
whose entries are floating point numbers, specially in the presence of round-off errors.
This fact was clearly stated in a historical paper by Golub and
Kahan \citep{GK}, where we can read the following sentences: ``In
the past the conventional way to determine the rank of $A$ was to
convert $A$ to a row-echelon form... But in floating-point
calculations it may not be so easy to decide whether some number
is effectively zero or not... {\it In other words, without looking
explicitly at the singular values there seems to be no
satisfactory way to assign rank to A}''.

Good general references in connection with the computation of the
SVD of a matrix are \citep{GVL} and \citep{D}.

\section{The inversion problem}
\label{}

Let $P(t)=(x(t), y(t))$ be a proper parametrization of a plane rational Bézier curve $C$  given by:
$$
x(t)=\frac{w_0a_0 \beta_0^{(n)}(t)+w_1a_1 \beta_1^{(n)}(t)+ \cdots + w_na_n \beta_n^{(n)}(t)}{w_0 \beta_0^{(n)}(t)+w_1 \beta_1^{(n)}(t)+ \cdots + w_n \beta_n^{(n)}(t)},
$$
$$
y(t)=\frac{w_0b_0 \beta_0^{(n)}(t)+w_1b_1 \beta_1^{(n)}(t)+ \cdots + w_nb_n \beta_n^{(n)}(t)}{w_0 \beta_0^{(n)}(t)+w_1 \beta_1^{(n)}(t)+ \cdots + w_n \beta_n^{(n)}(t)},
$$
$t \in [0,1]$.

\medskip
A parametrization $P(t)=(x(t), y(t))$ of a curve $C$ is said to be {\it proper} if every point on $C$ except a finite number of exceptional points is generated by exactly one value of the parameter $t$. Since every rational curve has a proper parametrization, we can assume that the considered parametrization is proper. Several results on the properness of curve parametrizations can be found in \citep{SW}.

\medskip
Our aim in this section is to solve the inversion problem for the case of a non exceptional given point $P_0=(x_0,y_0) \in C$, where $C$ is a Bézier curve properly parametrized by $P(t)$, that is to say, to compute the single value $t_0$ such that $P(t_0)=P_0$. The following result will be essential in our solution to this inversion problem.


\begin{thm}
Let $p(t)$ and $q(t)$ be two polynomials expressed in the Bernstein basis $\mathcal B_n$. If $t_0$ is a common root of $p(t)$ and $q(t)$, then the vector
 $(\beta_0^{(n-1)}(t_0), \beta_1^{(n-1)}(t_0), \ldots, \beta_{n-1}^{(n-1)}(t_0))$ is in the nullspace of the Bernstein-Bézout matrix of $p(t)$ and $q(t)$.
\end{thm}
\begin{pf}
Let $B$ be the Bernstein-Bézout matrix of $p(t)$ and $q(t)$, $\widehat B$ be the Bézout matrix of the polynomials $p(t)$ and $q(t)$ expressed in the monomial basis $\{1, t,t^2, \ldots, t^n  \},$ and $N$ be the lower triangular matrix of change of basis from the Bernstein basis $\mathcal B_{n-1}$ to the power basis $\{1, t,t^2, \ldots, t^{n-1}  \}$.
The following relationship is satisfied \citep{BG}:
$$
\widehat B=N B N^T.
$$
Let $t_0$ be a common root of $p(t)$ and $q(t)$, then
$$
\widehat B \left(
         \begin{array}{c}
           1 \\
           t_0 \\
           t_0^2 \\
           \vdots \\
           t_0^{n-1} \\
         \end{array}
       \right) = \left(
                   \begin{array}{c}
                     0 \\
                     0 \\
                     0 \\
                     \vdots \\
                     0 \\
                   \end{array}
                 \right)
$$
(see \citep{SZ}). In this way
$$
N B N^T \left(
         \begin{array}{c}
           1 \\
           t_0 \\
           t_0^2 \\
           \vdots \\
           t_0^{n-1} \\
         \end{array}
       \right) = \left(
                   \begin{array}{c}
                     0 \\
                     0 \\
                     0 \\
                     \vdots \\
                     0 \\
                   \end{array}
                 \right)
$$
and, as $N$ is a regular matrix
$$
B N^T \left(
         \begin{array}{c}
           1 \\
           t_0 \\
           t_0^2 \\
           \vdots \\
           t_0^{n-1} \\
         \end{array}
       \right) = \left(
                   \begin{array}{c}
                     0 \\
                     0 \\
                     0 \\
                     \vdots \\
                     0 \\
                   \end{array}
                 \right).
$$
Therefore the vector $$N^T \left(
         \begin{array}{c}
           1 \\
           t_0 \\
           t_0^2 \\
           \vdots \\
           t_0^{n-1} \\
         \end{array}
       \right) $$

belongs to the nullspace of the Bernstein-Bézout matrix $B$. Taking into account that $N$ is the matrix of change of basis from the Bernstein basis $\mathcal B_{n-1}$ to the power basis $\{1, t,t^2, \ldots, t^{n-1}  \}$

$$N^T \left(
         \begin{array}{c}
           1 \\
           t \\
           t^2 \\
           \vdots \\
           t^{n-1} \\
         \end{array}
       \right)=\left(
                 \begin{array}{c}
                   \beta_0^{(n-1)}(t) \\
                   \beta_1^{(n-1)}(t) \\
                   \beta_2^{(n-1)}(t) \\
                   \vdots \\
                   \beta_{n-1}^{(n-1)}(t) \\
                 \end{array}
               \right) $$
and in consequence the vector  $(\beta_0^{(n-1)}(t_0), \beta_1^{(n-1)}(t_0), \ldots, \beta_{n-1}^{(n-1)}(t_0))$ is in the nullspace of $B$, the Bernstein-Bézout matrix of $p(t)$ and $q(t).$
\end{pf}

\medskip

From now on, we will consider
$$
\begin{array}{l}
p(t) =w_0a_0 \beta_0^{(n)}(t)+w_1a_1 \beta_1^{(n)}(t)+ \cdots + w_na_n \beta_n^{(n)}(t) \\
~~~~~~~~ - x_0(w_0 \beta_0^{(n)}(t)+w_1 \beta_1^{(n)}(t)+ \cdots + w_n \beta_n^{(n)}(t)),\\
\\
q(t) =w_0b_0 \beta_0^{(n)}(t)+w_1b_1 \beta_1^{(n)}(t)+ \cdots + w_nb_n \beta_n^{(n)}(t) \\
~~~~~~~~ - y_0(w_0 \beta_0^{(n)}(t)+w_1 \beta_1^{(n)}(t)+ \cdots + w_n \beta_n^{(n)}(t)),
\end{array}
$$
and let $B$ be the Bernstein-Bézout matrix of these two polynomials $p(t)$ and $q(t)$.

Taking into account Theorem 2, the vector $(\beta_0^{(n-1)}(t_0), \beta_1^{(n-1)}(t_0),$ \break $\beta_2^{(n-1)}(t_0), \ldots, \beta_{n-1}^{(n-1)}(t_0))$ is in the nullspace of $B$. Moreover, as $P_0=(x_0,y_0)$ is a non exceptional point of $C$, which is properly parametrized by $P(t)$, the rank of $B$ is $n-1$ and therefore the vector $(\beta_0^{(n-1)}(t_0), \beta_1^{(n-1)}(t_0), \ldots, \beta_{n-1}^{(n-1)}(t_0))$ is a basis of the nullspace of $B$.

Let us observe that it is very usual not to know the point $P_0$ exactly. In this way, considering the results included in Section 2.2, it is convenient to compute an approximation of a basis of the nullspace of $B$ by using the SVD. Naturally, we do not have the exact matrix $B$ but an approximation of $B$ that we will denote by $\widetilde B$. A theoretical study of the application of SVD to the computation of an approximation of  the nullspace can be seen in Section 1 of Chapter 5 of \cite{ST}.

As the matrix $V$ of the SVD of $\widetilde B$ is orthogonal, the last column of $V$ gives us an approximation of a multiple of the vector $(\beta_0^{(n-1)}(t_0), \beta_1^{(n-1)}(t_0), \ldots, \beta_{n-1}^{(n-1)}(t_0))$ with Euclidean norm equal to 1:
$$
(z_0, z_1, \ldots, z_{n-1})=\alpha (\beta_0^{(n-1)}(t_0), \beta_1^{(n-1)}(t_0), \ldots, \beta_{n-1}^{(n-1)}(t_0)).
$$
Since
$$
\frac{z_i}{z_{i-1}}=\frac{\beta_i^{(n-1)}(t)}{\beta_{i-1}^{(n-1)}(t)} = \frac{{n-1 \choose i}(1-t)^{n-1-i}t^i}{{n-1 \choose i-1} (1-t)^{n-1-(i-1)}t^{i-1}}=\frac{(n-i)t}{i(1-t)}
$$
we obtain the value of the parameter $t$ we are interested in computing:
$$
t_0=\frac{i z_i}{i z_i+(n-i) z_{i-1}} \quad i=1,2, \ldots,n-1.
$$


\section{A more general situation}
\label{}

\bigskip
In this section we show how our procedure can also be used for the case of a non Bézier rational curve $C$ properly parametrized by $P(t)=(x(t),y(t))$, where $x(t)=\frac{u_1(t)}{u_2(t)}$, $y(t)=\frac{u_3(t)}{u_4(t)}$ and $u_i(t)$ ($i=1,\ldots,4$) are polynomials expressed in the Bernstein basis $\mathcal B_n$. In this situation some care must be taken when $deg(\bar p)\neq deg(\bar q)$, where  $\bar p(t)=u_1(t)-xu_2(t)$ and $\bar q(t)=u_3(t)-yu_4(t)$.

\medskip
Given a polynomial $p(t)$ expressed in the Bernstein basis $\mathcal B_n$, we denote by $deg(p)$ the {\it degree of $p(t)$ expressed in the power basis}. It must be noticed that $n\geq deg(p)$ (see Section 3 in \citep{BUG}).

\medskip
The following example \citep{MM07} serves to show what happens in this special situation:

\medskip
{\bf Example 1.} Let $\mathcal B_4$ be the Bernstein basis of the space of polynomials of
degree less than or equal to $4$, and let us consider the curve given by the parametric
equations
$$
x(t)= {4 \beta_0^{(4)}(t)+ 4 \beta_1^{(4)}(t)+ 3 \beta_2^{(4)}(t)+
3 \beta_3^{(4)}(t)+ 7 \beta_4^{(4)}(t)\over  \beta_0^{(4)}(t)+
\beta_1^{(4)}(t)+  \beta_2^{(4)}(t)+  \beta_3^{(4)}(t)+ 3
\beta_4^{(4)}(t)},
$$
$$
y(t)= 2 \beta_0^{(4)}(t)+ 3 \beta_1^{(4)}(t)+ 3 \beta_2^{(4)}(t)+
3 \beta_3^{(4)}(t)+ 4 \beta_4^{(4)}(t).
$$
In this case
$$
\bar p(t)=(4-x)\beta_0^{(4)}(t)+(4-x)\beta_1^{(4)}(t)+(3-x)\beta_2^{(4)}(t)+(3-x)\beta_3^{(4)}(t)+(7-3x)\beta_4^{(4)}(t)
$$
and
$$
\bar q(t)=(2 - y)\beta_0^{(4)}(t)+(3 - y)\beta_1^{(4)}(t)+(3 - y)\beta_2^{(4)}(t)+(3 - y)\beta_3^{(4)}(t)+(4 - y)\beta_4^{(4)}(t).
$$
However, if we write $\bar p$ and $\bar q$ in the power basis we have
$$
\bar p(t) = 4 - x - 6t^2 + 8t^3 + (-2x + 1)t^4
$$
(a polynomial of degree $4$ in $t$), while
$$
\bar q(t)= 2 - y + 4t - 6t^2 + 4t^3,
$$
a polynomial of degree $3$ in $t$.
In this example, our approach works perfectly when computing the parameter $t_0$ corresponding to every point $P_0=(x_0,y_0) \in C$ such that $x_0 \neq \frac{1}{2}$,  i.e., it works for every point $P_0 \in C$ except for the point $P_0=(1/2, -3.0395517)$. The reason why we have problems when $x_0=\frac{1}{2}$ is that $x_0=\frac{1}{2}$ makes $0$ the leading coefficient of $\bar p(t)$, and the nullspace of the corresponding Bernstein-Bézout matrix has dimension greater than $1$. When this happens, a  possible  solution is to consider the Sylvester matrix of the polynomials $\bar p(t)$ (expressed in the Bernstein basis $\mathcal B_4$) and $\bar q(t)$ (expressed in the Bernstein basis $\mathcal B_3$) \citep{WG} instead of the Bernstein-Bézout matrix of $\bar p(t)$ and $\bar q(t)$ expressed in the Bernstein basis $\mathcal B_4$, and then proceed in the same way as we have described for the Bernstein-Bézout matrix.


%

\section{Numerical examples}

In this section we include several examples illustrating the performance of our procedure. A detailed description of the way in which we proceed at each step of our approach is included in the first one.

\medskip
{\bf Example 2.} Let us consider a rational Bézier curve $C$ given by
$$
P(t)=(x(t), y(t))=\frac{\sum_{i=0}^{15}w_i (a_i,b_i){15 \choose i}t^i(1-t)^{n-i}}{\sum_{i=0}^{15}w_i {15 \choose i}t^i(1-t)^{n-i}},
$$
where the points of the control polygon are
{\small$$
\begin{array}{llll}
  (a_0,b_0)=(14,14), & (a_1,b_1)=(11,15), & (a_2,b_2)=(9,15), & (a_3,b_3)=(7,15), \\
  (a_4,b_4)=(4,14), & (a_5,b_5)=(3,12), & (a_6,b_6)=(3,10), & (a_7,b_7)=(7,8), \\
  (a_8,b_8)=(4,6), & (a_9,b_9)=(14,4), & (a_{10},b_{10})=(12,2), & (a_{11},b_{11})=(8,2), \\
  (a_{12},b_{12})=(6,2), & (a_{13},b_{13})=(4,3), & (a_{14},b_{14})=(3,4), & (a_{15},b_{15})=(2,5),
\end{array}
$$}
and the weights are
{\small$$
\begin{array}{cccccccc}
  w_0=2, & w_1=2, & w_2=2, & w_3=1, & w_4=2, & w_5=5, & w_6=5, & w_7=1, \\
  w_8=3, & w_9=3, & w_{10}=3, & w_{11}=3, & w_{12}=2, & w_{13}=1, & w_{14}=1, & w_{15}=1.
\end{array}
$$}
Let us consider $P_0=(8.50665, 14.3420)$, a point in $C$ that we know approximately (it is an approximation of the point corresponding to the value of the parameter $t=\frac{1}{7}$). We take $P_0=(x_0=\frac{850665}{10^5},y_0= \frac{143420}{10^4})$ and we construct the Bernstein-Bézout matrix of $p(t)$ and $q(t)$ by means of the {\it Maple} code introduced in Section 2.1 using exact arithmetic. The reason way we use exact arithmetic instead of floating point arithmetic is the following: although the algorithm presented in \citep{BG} for computing the entries of the Bernstein-Bézout matrix is fast (it has a computational cost of $O(n^2)$ arithmetic operations) this algorithm does not have high relative accuracy, and therefore it is possible that important errors appear in the computation of some entries of the matrix.

Once we have in {\it Maple} the approximation $\widetilde B$ of the Bernstein-Bézout matrix of $p(t)$ and $q(t)$, we put it in {\sc Matlab} where we compute its singular value decomposition. We include here the last column of the matrix $V$, which give us an approximation $\{(z_0, z_1, \ldots, z_{n-1})\}$ of the basis of the nullspace of $B$:
$$
V(:,15)=
\left(
\begin{array}{l}
-2.473682899590338e-001\\
   -5.771889609913881e-001\\
   -6.253061379846553e-001\\
   -4.169368183026464e-001\\
   -1.910461999315792e-001\\
   -6.368830764276734e-002\\
   -1.592182997982856e-002\\
   -3.033481337354997e-003\\
   -4.422976092589377e-004\\
   -4.912824640160427e-005\\
   -4.044767281991494e-006\\
   -2.444466877812876e-007\\
   -1.184192736731686e-008\\
   -5.593200217925554e-011\\
   -4.503186536883329e-011
\end{array}
  \right).
$$
In the computation of $t_0$ we use two consecutive components of the vector above chosen in the following way: we select the vector component which has the greatest absolute value, and then we take among the component before and after this the one that have the greatest absolute value. We proceed in this way because the components with greatest absolute value are usually the ones which are less sensitive to perturbations \citep{DDHK}. In this case we obtain
$$
t_0=\frac{2z_2}{2z_2+13z_1}=1.428606867264249e-001.
$$
When we compare the value of $t_0$ we have obtained with the exact value of the parameter $t=\frac{1}{7}$ of the exact point
$$
P=\Big(x\big(\frac{1}{7}\big), y\big(\frac{1}{7}\big)\Big)=\Big(\frac{78193109744768}{9191995131007}, \frac{131831466405881}{9191995131007} \Big) \in C,
$$
corresponding to the given point $P_0=(8.50665, 14.3420)$, we get that the relative error we made when using our procedure is
$$
\frac{\vert 1/7- t_0\vert}{1/7}= 2.48e-005.
$$
Taking into account that we have started from an approximation $P_0$ of the point $P\in C$ with only five exact digits, we can assert that the relative error that we have obtained by using our approach is small.

\medskip

Now we present a detailed description of the computations for the non Bézier curve presented in Example 1.

\medskip
{\bf Example 3.} Let us consider the curve $C$ introduced in Example 1 and a point $P_0=(x_0,y_0) \in C$ such that $x_0 \neq \frac{1}{2}$, for instance  $P_0=(3.5542169,2.8148148)$, which is an approximation of the point corresponding to the value of the parameter $t_0=1/3$. We show how our procedure works for these kind of points (all the points of $C$ except one).

If we consider $P_0=(\frac{35542169}{10^7}, \frac{28148148}{10^7})$, the singular values corresponding to the matrix $\widetilde B$ are:
$$\begin{array}{c}
  \sigma_1=4.212191730287018e+000, \\
  \sigma_2=2.075444341475023e+000, \\
  \sigma_3=5.981428444978487e-001, \\
  \sigma_4=3.357757839963324e-008.
\end{array}
$$
As only one singular value ($\sigma_4=3.357757839963324e-008$) can be considered $0$, the nullspace of $\widetilde B$ has dimension $1$ and therefore the last column of the matrix $V$ of the SVD of $\widetilde B$
$$
V(:,4)=\left(
         \begin{array}{c}
           5.111012703997072e-001 \\
            7.666518828223443e-001 \\
           3.833259078513024e-001\\
           6.388763832491218e-002\\
         \end{array}
       \right)
$$
gives us an approximation of the nullspace of the exact Bernstein-Bézout matrix $B$. Proceeding in the same way as in the previous example we obtain that
$$
t_0=\frac{z_1}{z_1+3z_0}=3.333333267311144e-001,
$$
a good approximation of $t=\frac{1}{3}$ if we take into account that $P_0$ is an approximation of the point corresponding to the value of the parameter $t_0=1/3$ with only $7$ exact digits.

Now we show what happens when we consider $P_0=(1/2, -3.0395517)$. In this case the singular values corresponding to the approximate Bernstein-Bézout  matrix $\widetilde B$  are
$$\begin{array}{c}
    \sigma_1=7.650996902942649e+001, \\
    \sigma_2=4.412324130498519e+001, \\
    \sigma_3=8.892019775556218e-009,\\
    \sigma_4=1.013508461498068e-015.
  \end{array}
$$
Since the singular values $\sigma_3$ and $\sigma_4$ can be considered $0$, the nullspace of $\widetilde B$ has dimension $2$, and therefore it is not possible to proceed by using the approach we have described in Section 3. As we have pointed out in Section 4, a possible solution for this situation could be to consider the Sylvester matrix of the polynomials $\bar p(t)$ (expressed in the Bernstein basis $\mathcal B_4$) and $\bar q(t)$ (expressed in the Bernstein basis $\mathcal B_2$) \citep{WG} instead of the Bernstein-Bézout matrix of $\bar p(t)$ and $\bar q(t)$ expressed in the Bernstein basis $\mathcal B_4$, and then proceed in the same way as we have described for the Bernstein-Bézout matrix.

\medskip
Finally we present an example with a Bézier curve of small degree.

\medskip
{\bf Example 4.} Let us consider the rational Bézier curve $C$ of degree $3$ given by
$$
P(t)=(x(t), y(t))=\frac{\sum_{i=0}^{3}w_i (a_i,b_i){3 \choose i}t^i(1-t)^{n-i}}{\sum_{i=0}^{3}w_i {3 \choose i}t^i(1-t)^{n-i}},
$$
where the points of the control polygon are
$$
\begin{array}{llll}
(a_0, b_0)=(1,9), & (a_1, b_1)=(2,1), & (a_2, b_2)=(5,1), & (a_3, b_3)=(4,1),
\end{array}
$$
and the weights are
$$
\begin{array}{llll}
w_0=1, & w_1=2, & w_2=2, & w_3=1.
\end{array}
$$
In this case, as the considered curve has small degree, it is easy to obtain (by using the {\it Maple} command {\tt ffgausselim}) the following explicit expression for an inversion formula:
$$
f(x,y)=\frac{81408x+13824yx-7680x^2-276552-1008y+1080y^2}{-216432-21024y+67200x+12672yx-6144x^2-2160y^2}.
$$
The exact point in $C$ corresponding to the parameter $t=\frac{1}{3}$ is $P=\big(\frac{8}{3}, \frac{109}{45}\big)$. When we consider its approximation $P_0=(2.66667, 2.42222)$ the value of $t_0$ that we get by means of our approach is $$3.333339104290224e-001,$$ while the value of $t_0$ we obtain by using in {\it Maple} the above inversion formula $f(x,y)$  for $P_0$ is $$3.333319852915495e-001.$$

These results show that, even in the cases we have an explicit inversion formula, the computation of the parameter $t_0$ by evaluating this function does not necessarily give better results than using the SVD.

\section{Conclusions}

In this work we have shown how the combination of classical algebraic geometry tools (resultant matrices) and numerical linear algebra tools (the singular value decomposition) is adequate for addressing the inversion problem for a rational Bézier curve, although the more general situation of rational non Bézier curves defined by polynomials in the Bernstein basis is also considered. The resultant matrices being used avoid the conversion between the Bernstein basis and the monomial basis.

It must be observed that the aim of our work was not the computation of an inversion formula, which is not an easy task for curves of high degree and must be approached by using different techniques. Instead we have considered the numerical problem of finding the parameter value corresponding to a given point on the curve, a point that usually in practice is not known exactly.

\section*{Acknowledgements}

This research has been partially supported by Research Grant
MTM2009-07315 from the Spanish Government.

\bibliographystyle{elsarticle-harv}


\begin{thebibliography}{00}

\bibitem[Bini and Gemignani, 2004]{BG} Bini, D. A., Gemignani, L., 2004. Bernstein-Bezoutian matrices. Theoretical Computer Science 315, 319--333.

\bibitem[Busé and Goldman, 2008]{BUG} Busé, L., Goldman, R., 2008. Division algorithms for Bernstein polynomials. Computer Aided Geometric Design 25, 850--865.

\bibitem[Busé and D'Andrea, 2006]{BDA} Busé, L., D'Andrea, C., 2006. A matrix-based approach to properness and inversion problems for rational surfaces. Applicable Algebra in Engenieering, Comunication and Computing 17, 393--407.

\bibitem[Chionh et al., 2002]{CZG} Chionh, E. W., Zhang, M., Goldman, R. N., 2002. Fast computation of the Bézout and Dixon resultant matrices. Journal of Symbolic Computation 33, 13-29.

\bibitem[Demmel, 1997]{D} Demmel, J. W., 1997. Applied Numerical Linear Algebra. SIAM.

\bibitem[Demmel et al., 2008]{DDHK} Demmel, J. W., Dumitriu, I., Holtz, O., Koev, P., 2008. Accurate and efficient expression evaluation and linear algebra. Acta Numerica 17, 87--145.

\bibitem[Farin, 2002]{FA} Farin, G., 2002. Curves and Surfaces for CAGD: A practical guide, fifth edition. Academic Press.

\bibitem[Farouki, 1991]{F} Farouki, R. T., 1991. On the stability of transformations between power and Bernstein polynomial forms. Computer Aided Geometric Design 8, 29--36.

\bibitem[Goldman et al., 1984]{GSA} Goldman, R. N., Sederberg, T. W., Anderson D. C., 1984. Vector elimination: A technique for the implicitization, inversion and intersection of planar parametric rational polynomial curves. Computer Aided Geometric Design 1, 327--356.

\bibitem[Golub and Kahan, 1965]{GK} Golub, G. H., Kahan, W., 1965. Calculating the singular values and pseudo-inverse of a matrix. SIAM Journal on Numerical Analysis, Ser. B, 2, 205--224.

\bibitem[Golub and Van Loan, 1996]{GVL} Golub, G. H., Van Loan, C. F., 1996. Matrix Computations, 3rd edition, Johns Hopkins University Press.

\bibitem[González-Vega and Rúa, 2009]{GVR} González-Vega, L., Rúa, I. F., 2009. Solving the implicitization, inversion and reparametrization problems for rational curves through subresultants. Computer Aided Geometric Design 26, 941--961.

\bibitem[Hoscheck and Lasser, 1993]{HL} Hoscheck, J., Lasser, D., 1993. Fundamentals of Computer Aided Geometric Design, A. K. Peters.

\bibitem[Marco and Martínez, 2004] {MM04} Marco, A., Martínez, J. J., 2004. A new source of structured singular value decomposition problems. Electronic Transactions on Numerical Analysis 18, 188--197.

\bibitem[Marco and Martínez, 2007] {MM07} Marco, A., Martínez, J. J., 2007. Bernstein-Bezoutian matrices and curve implicitization. Theoretical Computer Science 377, 65--72.

\bibitem[Sederberg et al., 1984]{SAG} Sederberg, T. W., Anderson, D. C., Goldman, R. N., 1984. Implicit representation of parametric curves and surfaces. Computer Vision, Graphics and Image Processing 28, 72--84.

\bibitem[Sederberg and Zheng, 2002]{SZ} Sederberg, T. W., Zheng, J., 2002. Chapter 15: Algebraic Methods for Computer Aided Geometric Design, in: Farin, G., Hoscheck J., Kim,  M. S. (Eds.), Handbook of Computer Aided Geometric Design, North-Holland.

\bibitem[Sendra and Winkler, 2001]{SW} Sendra, J. R., Winkler, F., 2001. Tracing index of rational curve parametrizations. Computer Aided Geometric Design 18, 771--795.

\bibitem[Strang, 1988]{ST} Strang, G., 1988. Linear Algebra and Its Applications, 3rd edition. Harcourt Brace Jovanovich.

\bibitem[Stewart, 1998]{STW} Stewart, G. W., 1998. Matrix Algorithms. Volume I: Basis Decompositions. SIAM.

\bibitem[Wang and Joe, 1995]{WJ} Wang, W., Joe, B., 1995. The geometric interpretation of inversion formulae for rational plane curves. Computer Aided Geometric Design 12, 469--489.

\bibitem[Winkler and Goldman, 2003]{WG} Winkler, J. R., Goldman R. N., 2003. The Sylvester resultant matrix for Bernstein polynomials, in: Lyche, T., Mazure, M. L., Schumaker, L. L. (Eds.), Curve and Surface Design: Saint-Malo 2002, 407--416. Nashboro Press.




\end{thebibliography}



\end{document}